\documentstyle[12pt,amsfonts,amscd]{amsart}
\newtheorem{Theorem}{Theorem}[section]
\newtheorem{Definition}[Theorem]{Definition}
\newtheorem{Proposition}[Theorem]{Proposition}

\newtheorem{Lemma}[Theorem]{Lemma}
\newtheorem{Corollary}[Theorem]{Corollary}
\theoremstyle{remark}
\newtheorem{Example}[Theorem]{Example}

\def\Aut{\operatorname{Aut}}
\def\dim{\operatorname{dim}}

\def\be{\begin{enumerate}}
\def\ee{\end{enumerate}}
\def\bT{\begin{Theorem}}
\def\eT{\end{Theorem}}
\def\bP{\begin{Proposition}}
\def\eP{\end{Proposition}}
\def\bD{\begin{Definition}}
\def\eD{\end{Definition}}
\def\bE{\begin{Example}}
\def\eE{\end{Example}}
\def\bL{\begin{Lemma}}
\def\eL{\end{Lemma}}
\def\bC{\begin{Corollary}}
\def\eC{\end{Corollary}}

\begin{document}
\title{Fixed points and Determining Sets for Holomorphic Self-Maps of a Hyperbolic Manifold}
\author{Buma L. Fridman, Daowei Ma and Jean-Pierre Vigu\'{e}}
\begin{abstract} We study fixed point sets for holomorphic automorphisms (and endomorphisms) on complex manifolds. The main object of our interest is to determine the number and configuration of fixed points that forces an automorphism (endomorphism) to be the identity. These questions have been examined in a number of papers for a  bounded domain in ${\Bbb C}^n$. Here we resolve the case for a general finite dimensional hyperbolic manifold. We also show that the results for non-hyperbolic manifolds are notably different.
\end{abstract}
\address{ fridman@@math.wichita.edu, Department of Mathematics,
Wichita State University, Wichita, KS 67260-0033, USA}
\address{ dma@@math.wichita.edu, Department of Mathematics,
Wichita State University, Wichita, KS 67260-0033, USA}
\address{ vigue@@math.univ-poitiers.fr, UMR CNRS 6086, Universit\'{e} de Poitiers, 
Math\'{e}matiques, SP2MI, BP 30179, 86962 FUTUROSCOPE, FRANCE}
\maketitle \setcounter{section}{-1}
\section{Introduction}

Let $M$ be a complex manifold. $H(M,M)$ is the set of holomorphic maps from $M$ to $M$, i.e., the set of endomorphisms of $M$. A special case of endomorphisms are automorphisms of $M$, $Aut(M)\subset H(M,M)$.

\begin{Definition}
A set $K\subset M$ is called a determining subset of $M$ with respect to $
Aut(D)$ ($H(M,M)$ resp.) if, whenever $g$ is an automorphism (endomorphism resp.) such that $g(k)=k$ $\forall k\in K$, then $g$ is the identity map of $M$.
\end{Definition}

The notion of a determining set was first introduced in a paper written by the first two authors in collaboration with Steven G. Krantz and Kang-Tae Kim [FK1]. That paper
was an attempt to find a higher dimensional analog of the following result
of classical function theory [PL]: if $f:M\rightarrow M$ is a conformal
self-mapping of a plane domain $M$ which fixes three distinct points then $
f(\zeta )=\zeta $. \newline
This one-dimensional result is true even for endomorphisms of a bounded domain $D\subset\subset \Bbb C$. To prove this one needs to first use the well known theorem, stating that if an endomorphism of $D$ fixes two distinct points, then it is an automorphism; and then use the above [PL] theorem.

Determining sets (for automorphisms and endomorphisms) in case of bounded domains in ${\Bbb C}^n$ have been further investigated in the following papers [FK2], [KK], [Vi1], [Vi2], [FM].

Let $W_s(M)$ denote the set of $s$-tuples $(x_1,\dots, x_s)$, where $x_j\in M$, such that $\{x_1,\dots, x_s\}$ is a determining set with respect to $Aut(M)$. Similarly, $\widehat{W}_s(M)$ denotes the set of $s$-tuples $(x_1,\dots, x_s)$ such that $\{x_1,\dots, x_s\}$ is a determining set with respect to $H(M,M)$.
So $\widehat{W}_s(M)\subseteq W_s(M)\subseteq M^s$.
We now introduce two values $s_0(M)$ and $\widehat{s}_0(M)$. In case $Aut(M)={id
}$, $s_0(M)=0$, otherwise $s_0(M)$ is the least integer $s$, such that $W_s(M)\neq \emptyset $. If $W_s(M)=\emptyset $ for all $s$ then $s_0(M)=\infty $. Analogously symbol $\widehat{s}_0(M)$ denotes the least integer $s$ such that $\widehat{W}_s(M)\neq \emptyset $, if no such integer exists (i.e. $\widehat{W}_s(M)=\emptyset $ for all $s$) then $\widehat{s}_0(M)=\infty $. In all cases $s_0(M)\leq \widehat{s}_0(M)$.
\newline The main objectives of this paper are: first, to generalize the results
for bounded domains in ${\Bbb C}^n$ to hyperbolic manifolds, and second, to illustrate that for the non-hyperbolic manifolds the results are quite different. 
\newline The Bergman metric on a bounded domain in ${\Bbb C}^n$ proved quite useful for the investigation related to determining sets. Such a Riemannian metric however is not always available on a hyperbolic manifold $M$; to overcome this obstacle we construct for any point $x\in M$ an invariant (with respect to $Aut(M)$) Hermitian metric in a neighborhood (open but not necessarily connected) of that point. 
\newline The paper is roughly divided into three parts. First we introduce the
Hermitian metric mentioned above. Second, we completely resolve the case for
a hyperbolic manifold. Third, we prove (two) theorems to show that the
case of non-hyperbolic manifolds is remarkably different.
\newline Here is a brief description of the second and third parts of the paper. 
In [Vi2] the estimate $\widehat{s}_0(D)\leq n+1$ was established  for all bounded domains in ${\Bbb C}^n$. In Section 2 we generalize this by proving the same inequality for hyperbolic manifolds of dimension $n$. This certainly implies same inequality for {\it automorphisms} of a hyperbolic manifold $M$, $s_0(M)\leq n+1$. However for automorphisms much more information can be provided. $s_0(M)$ depends on how large the group $Aut(M)$ is, and corresponding estimates on $s_0(M)$ are given in section 3.
In section 4 we show that if $dim(M)=n$ then the general estimate ($s_0(M)\leq n+1$) can be refined to $s_0(M)\leq n$ for domains that are not biholomorphic to the unit ball $B^n\subset {\Bbb C}^n$ (i.e. the only hyperbolic manifolds for which $s_0(M)=n+1$ are those biholomorphic to the ball). 
\newline If a positive integer $s\geq s_0(M)$, then $W_s(M)\neq \emptyset $, so there
are $s$ points such that if an automorphism of $M$ fixes these points it
will fix any point of $M$. Now the question arises whether the choice of
these $s$ points is generic. The answer is positive for any hyperbolic manifold $M$:  $W_s(M)\subseteq M^s$ is open and dense if not empty(section 5).
\newline Similar topological properties for the determining sets of {\it endomorphisms} of a general hyperbolic manifold do not hold. We address related questions in the concluding part of section 5.
\newline Section 6 is devoted to examining the situation for {\it non-hyperbolic} manifolds. We first give a complete description of the value for $\widehat{s}_0(M)$ for a one-dimensional manifold $M$ (theorem \ref{od}). Then for higher dimensional manifolds we prove that $\widehat{s}_0(M)=\infty $ for a general Stein
manifold $M$ that has the following property: any finite number of points
lie in a one dimensional submanifold (for the precise statement see theorem \ref{sm} )
\section{Construction of a locally invariant Hermitian metric}
Our main effort in this section will be the construction of a locally invariant (with respect to the automorphism group) metric in a neighborhood of any point in a general hyperbolic manifold. First we present some preliminary statements.
\newline Throughout this section $M$ denotes a hyperbolic manifold of finite dimension, $Aut(M)$ is its group of holomorphic automorphisms.
 
\bL\label{nor} $Aut(M)$ is a normal family.
\eL
Various versions of this statement have been used before. However, we cannot find a direct reference to this result in the literature. Therefore a brief proof is presented here.
\begin{pf}

It suffices to prove that if $x_0\in M$, if $f_j\in Aut(M)$ is a sequence such that the closure $Q$ of the set $\{f_j(x_0): j\in {\Bbb N}\}$ is compact, and if $K$ is a compact subset of $M$, then
$$S:=\cup_{j=1}^{\infty} f_j(K)\subset\subset M.$$
Let $d(\cdot,\cdot)$ denote the Kobayashi distance. For $x\in M, r>0$ let $b(x,r)=\{y\in M: d(x,y)<r\}$. Let $\psi(x)=\sup\{r>0: b(x,r)\subset\subset M\}$. Now we set
$$m=\max\{d(x_0,x): x\in K\},\;\;\; \delta=\min\{\psi(x): x\in K\},$$
and 
$$ P=\{x\in M: d(x, Q)\le m, \psi(x)\ge \delta\}.$$
Then $P$ is compact and $S\subset P$.
\end{pf}

Now we note the following. Let $a\in M$, $f:M\rightarrow M$ a holomorphic map such that $f(a)=a$.  Consider a small Kobayashi ball $b=b(a,\varepsilon )\,$that is biholomorphic to a bounded domain in  ${\Bbb C}^n$, and whose closure is compact in $M$. Since the Kobayashi distance is non-increasing under holomorphic maps, we have $f:b\rightarrow b$. If $f\in Aut(M)$, then $f\mid _b\in Aut(b)$. The following three statements (cf. [Vi1]) hold for bounded domains in ${\Bbb C}^n$; by using this remark one can prove them for any hyperbolic manifold.
\bL \it{Let $a\in M$, $f:M\rightarrow M$ a holomorphic map such that $f(a)=a$ and $f^{\prime }(a)=id$. Then $f=id$.}
\eL

\bL Let $a\in M$, $f\in Aut(M)$ and  $f(a)=a$. Then all the eigenvalues of  $
f^{\prime }(a)$ are of modulus one, and the matrix $f^{\prime }(a)$ is
diagonalizable.
\eL

\bC\label{ei}In the assumption of the above Lemma, if $f \ne id$, one can find an appropriate power $k$ such that the $k$-th iteration of $f$, $f^k=h\in Aut(M)$ will have the following properties: $
h(a)=a,h^{\prime }(a)$ has at least one eigenvalue with non-positive real part. 
\eC

Let $z\in M$. Below we use the notion of an isotropy group $I_z(M)= \{g\in Aut(M):g(z)=z\}$.
\bL[H. Cartan]([Ca1, p.80])
{Let $D\subset \subset {\Bbb C}^n$, let $z\in D$, and let $I_z=I_z(D)$
be the isotropy subgroup at $z$ of the automorphism group of $D$. Then there
exists a holomorphic map $\phi :D\rightarrow {\Bbb C}^n$ such that $\phi
(z)=0,$  $\phi ^{\prime }(z)=id$, and for all $f\in I_z$
one has $\phi \circ f=f^{\prime }(z)\circ \phi $.}
\eL

\noindent  As in [Vi1, thm 2.3], for the proof of this Lemma, we define $\phi :D\to 
{\Bbb C}^n$ by 
\begin{equation*}
\phi (\zeta )=\int_{G_z}f^{\prime }(z)^{-1}(f(\zeta )-z)\,d\mu (f),
\end{equation*} 
where $d\mu $ is the Haar measure on $I_z$. Then $\phi (z)=0$, $\phi ^{\prime }(z)=id$ (and therefore $\phi$ is locally biholomorphic),
and $\phi \circ g=g^{\prime }(z)\circ \phi $ for each $g\in I_z$.
\newline
\newline Let $M$ again be a hyperbolic manifold, $x\in M$, $T_xM$ the tangent space of $M$
at $x$, $I_x=I_x(M)$ is the isotropy subgroup fixing $x$. The compact
group\thinspace $I_x$ acts on $T$ as differential maps: for $g\in I_x$ , $
v\in T$ , $g_*(v)=dg(x)v$. Since the above Lemma can be considered in a small
neighborhood of $x$, and $T$ is isomorphic to ${\Bbb C}^n$ the following statement holds.
\bL\label{ca} For any point $x\in M$ there exists a small neighborhood $V\ni x$, such that there is an injective holomorphic map $\phi: V\rightarrow T$ such that $g_*\circ \phi=\phi\circ g$ for $g\in I_x$, and $d\phi(x)=id$, the identity map of $T=T_xM$. 
\eL

Finally we introduce an Hermitian invariant metric on a neighborhood of any point in $M$.
\bL\label{me} Let $M$ be a hyperbolic manifold, let $G=Aut(M)$,  and let $x\in M$. Then there is a neighborhood $U$ of $x$ such that $G(U)=U$, and a $C^\infty$ Hermitian metric on $U$ that is invariant under $G$.
\eL

\begin{pf} Since $M$ is hyperbolic, the automorphism group $G$ is a Lie group (see [Ko]) and the isotropy group $I_x$ is a compact subgroup of $G$. The orbit $G(x)$ is an embedded submanifold of $M$. 
Let $T=T_xM$ be the tangent space of $M$ at $x$. Then $T$ is a complex vector space and is isomorphic to ${\Bbb C}^n$. The elements of the compact group $I_x$ act on $T$ as differential maps: for $g\in I_x$, $g_*(v)=dg(x)v$. Let $h$ be a Hermitian metric on $T$ invariant under $I_x$. By Lemma \ref{ca}, there exist a small neighborhood $V$ of $x$ in $M$ and an injective holomorphic map $\phi: V\rightarrow T$ such that $g_*\circ \phi=\phi\circ g$ for $g\in I_x$, and $d\phi(x)=id$, the identity map of $T=T_xM$. The real subspace $P$ of $T$ consisting of vectors tangent to $G(x)$ is invariant under $I_x$. So the orthogonal complement (with respect to the real part of $h$) $Q$ of $P$ is also invariant under $I_x$.  Let $S_1=\{v\in Q: \|v\|<\delta\}$, where $\|\cdot\|$ is the norm induced by the Hermitian metric $h$, and choose $\delta>0$ so small that $S_1\subset\subset \phi(V)$. Note that $S_1$ is invariant under $I_x$. Let $S=\phi^{-1}(S_1)$. Then $I_x(S)=S$. Furthermore, for $g\in G$, $g(S)\cap S\not= \emptyset$ iff $g\in I_x$. The tube $G(S)$ is diffeomorphic to the the normal bundle of $G(x)$ in $M$ and to the twisted product $G\times_{I_s} S$. The pull-back $h_0=(\phi|_S)^{\ast} h$ is a Hermitian metric on the restriction to $S$ of the tangent bundle $TM$. Now we define a Hermitian metric $h_1$ on $U=G(S)$ as follows. If $y\in U$ and $u, v\in T_y$, then there is a $g\in G$ such that $g(y)\in S$, and we define $h_1(u, v)=h_0(g_{\ast}u, g_{\ast}v)$. One can see that $h_1$ is well-defined, since if $g(y), g'(y)\in S$, then $g' g^{-1}\in I_x$. Now $h_1$ is a $C^\infty$ metric on $U$ that is invariant under $G$. 
\end{pf}

\section{An estimate for $\widehat{s}_0(M)$}
We need the following lemma (Thm. 5.2 in [Vi2])
\bL Let $D$ be a bounded domain in ${\Bbb C}^n$, $a\in D$. Then there is an open 
$U\subset D^n$ such that $(a,...,a)\in \overline{U}$ and for all $(z_1,...,z_n)\in U$, $(a,z_1,...,z_n)\in \widehat{W}_{n+1}(D)$.
\eL
\bT\label{mes} Let $M$ be a hyperbolic manifold of complex dimension $n$. Then $\widehat{s}_0(M)\le n+1$.
\eT
\begin{pf}
Pick a point $a\in M$. Let $f:M\rightarrow M$ be a holomorphic map such that $f(a)=a$. Consider a small Kobayashi ball $b=b(a,\varepsilon )\,$ whose closure is compact in $M$, and such that $b$ is biholomorphic to a bounded domain $D$ in  ${\Bbb C}^n$; let $h:b\rightarrow D$ be such a biholomorphic map. Note that since the Kobayashi distance is non-increasing under holomorphic maps, we have $f:b\rightarrow b$, and therefore $g=h\circ f\circ h^{-1}:D\rightarrow D$. By using the preceding lemma, one
can pick n points $z_1,...,z_n\in D$, such that $Z=(h(a),z_1,...,z_n)\in 
\widehat{W}_{n+1}(D)$.\thinspace Consider the set of $n+1$ points $%
h^{-1}(Z)=(a,h^{-1}(z_1),...,h^{-1}(z_n))\subset b$. If our function $f\in
H(M,M)$ (in addition to $a$) is also fixing all points $h^{-1}(z_j)$, i.e.  $%
f\mid _{h^{-1}(Z)}=id$, then $g\mid _Z=id$ and therefore $g=id$. We conclude
that  $f\mid _b=id,$ and consequently $f=id$. So, $h^{-1}(Z)\in \widehat{W}%
_{n+1}(M)$, and therefore $\widehat{s}_0(M)\le n+1$.
\end{pf}

\section{Estimates for $s_0(M)$}
The goal of this section is to provide estimates for $s_0(M)$ for a hyperbolic manifold $M$, $dim(M)=n$.
\newline Since $s_0(M)\leq \widehat{s}_0(M)$ theorem \ref{mes} implies
\newline \newline {\it For any hyperbolic manifold $M$ of complex dimension $n$, $s_0(M)\le  n+1$}.
\newline \newline
{\it Remark.} In the next section we prove a refined inequality $s_0(M)\le  n$ for $M$ not biholomorphic to the unit ball in ${\Bbb C}^n$.
\newline \newline
If $H$ is (isomorphic to) a subgroup of the unitary group $U(n)$, let $k(H)$ denote the least number $k$ of vectors $u_1,\dots, u_k$ such that if $h\in H$ and if $h(u_j)=u_j$ for $j=1,\dots,k$ then $h=id$. For $x\in M$ the isotropy group $I_x(M)$ is isomorphic to the group of its differentials at $x$, and these differentials are unitary with respect to the locally defined Hermitian inner product (the existence of which was proved in Lemma \ref{me}) on the tangent space $T_x(M)$. So $I_x(M)$ is isomorphic to a subgroup of $U(n)$.

\bT\label{li}
$s_0(M)\le 1+\min\{k(I_x(M)): x\in M\}$.
\eT

\begin{pf} Choose $x\in M$ so that $k(I_x(M))=\min\{k(I_x(M)): x\in M\}$. Denote that number by $k$. Let $u_1,\dots, u_k$ be vectors in $T_x M$ such that if $h\in I_x(M)$ and if $dh(x)(u_j)=u_j$ for $j=1,\dots,k$ then $dh=id$ (hence $h=id$). For each $u_j$, let $x_j$ be a point on the geodesic through $x$ in the direction $u_j$, so close to $x$ that the geodesic is the unique length minimizing geodesic from $x$ to $x_j$. Let $f$ be an automorphism of $M$ fixing $x, x_1, \dots, x_k$. Then $df(x)$ fixes $u_1, \dots, u_k$. It follows that $df(z)=id$ and $f=id$. Therefore, $s_0(M)\le  1+\min\{k(I_x(M)): x\in M\}$.
\end{pf}

Let $G$ be a subgroup of $Aut(M)$. By $s_0(M,G)$ we denote the minimum
number of distinct points in $M$ such that if $g\in G$, and $g$ fixes all
these points, then $g=id$. So, $s_0(M)=s_0(M,Aut(M)).$

\bL\label{es} Let $M$ be a hyperbolic manifold, let $G$ be a subgroup of $Aut(M)$, and let $q=\dim G$. If $q\ge 1$, then $s_0(M, G)\le q$. If $q=0$, then $s_0(M, G)\le 1$.
\eL

\begin{pf} First we consider the case where $q\le 1$. Let $e$ denote the identity element of $G$, and let $Q=G\backslash \{e\}$. For each $g\in Q$, the set $\{x\in M: g(x)=x\}$ is an analytic set of $M$ of dimension $\le 2n-2$. The set $W_1:=\{(g,x)\in Q\times M: g(x)=x\}$ is an analytic set of $Q\times M$ of dimension $\le (2n-2)+q\le 2n-1< dim M$. Let $W$ denote the set of fixed points of nontrivial elements of $G$. Since $W=\pi(W_1)$, where $\pi: Q\times M\to M$ is the projection, and since $dim W_1< dim M$, we see that $W\ne M$. Therefore, $s_0(M, G)\le 1$. 

Now we assume that $q\ge 2$. There must be an orbit $Q$ of $G$ of positive dimension. Let $x\in Q$, and let $H:=G_x$ be the subgroup of $G$ consisting of elements $g$ satisfying $g(x)=x$. Then $\dim H< \dim G$. By induction hypothesis, $s_0(M, H)\le \dim G-1$. Therefore, $s_0(M, G)\le 1+s_0(M, H)\le \dim G$.
\end{pf}
As a corollary we get
\bT  If 
$\dim (Aut(M))\ge 1$, then $s_0(M)\le \dim (Aut(M))$. If $\dim (Aut(M))=0$, then $s_0(M)\le 1$.
\eT

\section{A characterization of the ball in ${\Bbb C}^n$}
This section is devoted to the proof of the following statement (which is a generalization of Theorem 1.1 in [FM]).
\bT\label{bt} Let $M$ be a hyperbolic manifold of dimension $n$. $s_0(M)=n+1$ if and only if $M$ is biholomorphic to the unit ball $B^n$ in ${\Bbb C}^n$.
\eT
The estimate  $s_0(B^n)=n+1$ can be easily verified (see for example [FM]).

The rest of this section will be devoted to the proof that $s_0(M)=n+1$ implies that $M$ is biholomorphic to the unit ball. 
To prove this we need the following two lemmas.

\bL\label{tr} Let $M$ be a hyperbolic manifold and $x\in M$. Suppose that the isotropy group $I_x$ is transitive on the (real) directions at $x$. Then $M$ is biholomorphic to the unit ball in ${\Bbb C}^n$.
\eL

\begin{pf} Since $I_x$ is transitive on the directions at $x$, the group $Aut(M)$ is not finite. Since the automorphism group of a compact hyperbolic manifold must be finite (see [Ko, p.~70]), we see that $M$ is noncompact. By the main theorem in [GK], $M$ is biholomorphic to ${\Bbb C}^n$.
\end{pf}

For a subgroup $H$ of the unitary group $U(n)$ we use the notion $k(H)$ introduced at the beginning of section 3. The following Lemma was proved in ([FM], Lemma 1.4).

\bL\label{tri}
If $H$ is a subgroup of $U(n)$ with $n\ge 2$ and if $H$ is not transitive on $S^{2n-1}$ then $k(H)\le n-1$.
\eL

We are now ready to prove the remaining portion of Theorem~$\ref {bt}$ (i.e., $s_0(M)=n+1$ implies that $M$ is biholomorphic to the unit ball).

\begin{pf} So, let $s_0(M)=n+1$. If $n=1$ the statement ($M$ is biholomorphic to the unit disc $B^1$) is true. Indeed, if $M$ is not biholomorphic to the disc or the annulus, its automorphism group is discrete. For each element $g \in Aut(M)$, $g\ne id$ the set of fixed points is discrete. Therefore there is a point $x\in M$ that is not a fixed point of any nontrivial automorphism. This point will then form a determining set, and so, $s_0(M)\le 1$. For the annulus $s_0(M)=1$. Therefore if $s_0(M)=2$, $M$ is biholomorphic to the unit disc.

Consider now the case where $n\ge 2$. Let $z\in M$. Suppose that $M$ is not biholomorphic to $B^n$. Then $I_z(M)$ is not transitive on the directions at $z$, by Lemma $\ref {tr}$. Since $I_z(M)$ is (isomorphic to) a subgroup of $U(n)$, by Lemma $\ref {tri}$, $k(I_z(M))\le n-1$. It follows (see Theorem $\ref {li}$) that $s_0(M)\le 1+ k(I_z(M))\le n$ if $M$ is not biholomorphic to $B^n$.
\end{pf}

\section{Determining sets $W_s(M)$ are open and dense}

Our aim in this section is to prove the following theorem.

\bT\label{op} Let $M$ be a hyperbolic manifold and $s\ge 1$. Then $W_s(M)\subset M^s$ is open;
if in addition $W_s(M)\neq \emptyset $, then $W_s(M)$ is dense in $M^s$.
\eT
Denote $W=W_s(M)$. First we prove that $W\subset M^s$ is open.
\begin{pf} Suppose $W$ is not open. Then one can find a sequence of $s$-tuples $Z_j=(x_1,...,x_s^j)\in M^s$
that converges to $Z=(x_1,...,x_s)\in M^s$ and such that $Z_j$ is not a
determining set for $M$, and $Z$ is. For each $j$ there is an $f_j\in Aut(M),
$  $f_j\mid _{Z_j}=id$, but $f_j\neq id$. By Corollary \ref{ei}
(replacing $f_j$ by an appropriate iteration of $f_j$ if needed) we may assume that the
real part of at least one eigenvalue of $f_j^{\prime }(x_1^j)$ is
non-positive. Switching again to a subsequence, if necessary, we find a
sequence of automorphisms whose limit (see Lemma \ref{nor}) is $g\in Aut(M)$, such
that $g\mid _Z=id$ , and one of the eigenvalues of  $g^{\prime }(x_1)$ is
non-positive. Therefore $g\neq id$ which contradicts the original assumption
that  $Z$ is a determining set for $M$. 
\end{pf}
{\it Remark}. The above proof of the theorem for a bounded domain is given in [Vi1, Thm 3.1]. One can also prove Theorem \ref{op} by using the idea of [FM, Lemma 2.3]. \newline
\newline Now suppose that $W\neq \emptyset $. We need to prove that $W$ is dense in $M^s$.
\newline First we introduce some notation. If $G$ is a subgroup of $\Aut(M)$, 
$W_s(M,G)$ denotes the set of $s$-tuples $(x_1,\dots, x_s)$, where $x_j\in M$, such that each element $g\in M$ satisfying $g(x_j)=x_j$ for $j=1,\dots,s$ has to be the identity.
\newline Let $\rho_x(\cdot, \cdot)$ denote the metric introduced in Lemma \ref{me} for a point $x\in M$. Let $b(x,r)$ denote the ball with center $x$ and radius $r$ in that metric. Let $\overline b(x,r)$ be the closure of $b(z,r)$ in $M$.

\bL Suppose that $G$ is a subgroup of $Aut(M)$. If $W_1(M,G)\neq \emptyset $ then $W_1(M,G)$ is dense in $M$.
\eL
\begin{pf} In this proof, let $W=W_1(M,G)$. Suppose that $W$ is not dense in $M$. Then the closure $K$ of $W$ in $M$ is not equal to $M$. Let $p$ be a boundary point of $K$ in $M$. Denote $\rho(\cdot, \cdot)=\rho_p(\cdot, \cdot)$. Choose $r>0$ such that the closure of $b(p, 4r)$ in $U$ is compact, where $U$ is a neighborhood from Lemma \ref{me} (chosen for the point $p$), and such that each pair of points of $b(p, 4r)$ is connected by a unique length-minimizing geodesic segment in that metric. There exist points $z, w$ such that $\rho(z, p)<r$, $\rho(w,p)<r$, $w\in W$, and $z\not\in K$. Note that the orbit of $w$, $G(w)\subset W$. Let $Q=G(w)\cap \overline b(p, 4r)$. Then $Q$ is compact and $Q \subset W$. Let $u$ be a point of $Q$ nearest to $z$. Then $u$ is also a point of $G(w)$ nearest to $z$, and $R:=\rho(z, u)\le \rho(z,w)<2r$. Choose a point $y$ on the unique length-minimizing geodesic segment from $z$ to $u$ such that $y\not\in K$ and $y\ne z$. For each point $x$ of $G(w)$, we see that 
$$\rho(z, y)+\rho(y, x)\ge \rho(z, x)\ge \rho(z, u),$$
and that the two equalities hold simultaneously only if $x=u$. Hence, $\rho(z, y)+\rho(y, x)> \rho(z, u)=R$ for each $x\in G(w)$, $x\ne u$. It follows that $\rho(y,x)> R-\rho(z,y)=\rho(y,u)$ for each $x\in G(w)$, $x\ne u$. Therefore, $u$ is the unique point of $G(w)$ nearest to $y$. Since $y\not\in K$, there is a nontrivial $g\in G$ such that $g(y)=y$. Now $\rho(y, u)=\rho(g(y), g(u))=\rho(y, g(u))$ forces $g(u)=u$. Since $u\in W$, the map $g$ must be the identity, contradicting the fact that $g$ is not trivial. Therefore, $W_1(M,G)$ is dense in $M$.
\end{pf}
{\it Proof} of Theorem \ref{op}.
We have already proved that $W_s(M)$ is open in $M^s$. Suppose now that $W_s(M)\neq \emptyset $.  For $g\in Aut(M)$ let $Q_s(g)$ denote the mapping 
\[
Q_s(g):M^s\rightarrow M^s,\;\;Q_s(g)(z_1,\dots ,z_s)=(g(z_1),\dots ,g(z_s)).
\]
Let $G=\{Q_s(g):g\in Aut(M)\}$. Then $G\subset Aut(M^s)$, and $
W_1(M^s,G)=W_s(M)$. By the previous lemma, $W_s(M)$ is dense in $
M^s$.
\newline \newline By using the same approach as in Theorem 5.1 in [Vi2] one can establish the following
\bT If  $M$ is a taut manifold then $\widehat{W}_s(M)$ is open in $M^s$ for all $s\ge 1$.
\eT
In general $\widehat{W}_s(M)$ does not have to be open in $M^s$ (see [FM]), nor be dense in $M^s$ (cf. [Vi2],[FM]).
\section{Results concerning the non-hyperbolic case}

\subsection{One dimensional manifolds} 

\bT\label{od} For a one dimensional complex manifold $M$, $\hat s_0(M)=2$ or $\infty$. More precisely, if $M$ is holomorphically equivalent to the complex plane $\Bbb C$, the truncated complex plane ${\Bbb C}^*$, or the Riemann sphere $\Bbb P$, then $\hat s_0(M)=\infty$; otherwise $\hat s_0(M)=2$. \eT

\begin{pf} It is well known that either $M$ is biholomorphic to ${\Bbb C}$, ${\Bbb C}^*$, $\Bbb P$, or a torus, or else it is a hyperbolic manifold. 

Suppose $S=\{x_1, \dots, x_k\}$ is a finite set in $\Bbb C$. Choose $y\in {\Bbb C}\setminus S$. Let $f$ be a polynomial such that $f(x_j)=x_j$ and $f(y)=y+1$. Then $f$ is a nonidentity holomorphic self map of ${\Bbb C}$ fixing each point of $S$. So $\hat s_0({\Bbb C})=\infty$.

We now consider $\Bbb P$. Note that the map $f$ in the last paragraph is also a holomorphic self map of ${\Bbb P}$ fixing the point at infinity. We see that $\hat s_0({\Bbb P})=\infty$.

Suppose $S=\{x_1, \dots, x_k\}$ is a finite set in ${\Bbb C}^*$. Choose $y_j$ so that $x_j=\exp y_j$. Let $g$ be a polynomial such that $g(x_j)=y_j$ and let $f(z)=\exp g(z)$. Then $f$ is a nonidentity holomorphic self map of ${\Bbb C}^*$ fixing each point of $S$. So $\hat s_0({\Bbb C}^*)=\infty$.

Consider a torus $T$ corresponding to a lattice $L$  in the complex plane. 
Let $\pi: {\Bbb C}\rightarrow T$ be the projection. It is well known that each holomorphic self map of $T$ has the form $f(\pi(z))=\pi(\lambda z+b)$, where $b\in {\Bbb C}$ and $\lambda\in \Lambda:=
\{x\in {\Bbb C}: xL\subset L\}$. Clearly $f$ is the identity iff $\lambda=1$ and $b\in L$. Let $F$ be the field generated by $L\cup \Lambda$. Then $F$ is countable. Choose $r\in {\Bbb C}\setminus F$. Let $x=\pi(0), y=\pi(r)$. Suppose that $f(x)=x$ and $f(y)=y$. Then
$$\lambda\cdot 0+b=0+ p,\;\;\; \lambda r+b=r+q,$$
for some $p, q\in L$. It follows that $b\in L$ and 
$$(\lambda-1)r=q-p.$$
Now $(\lambda-1)\in F$, $(q-p)\in F$, but $r\not\in F$. It follows that $\lambda-1=0$ and $f=id$. Therefore, $\hat s_0(T)=2$.

If $M$ is a hyperbolic manifold $M$ of dimension one, then $\hat s_0(M)=2$ by Theorem 2.2.
\end{pf}

\subsection{Higher dimensional manifolds}
The main statement in this section is the following
\bT\label{sm} Let $M$ be a Stein manifold,  $\dim (M)\geq 2$, and such that for any $k$ distinct points $\{x_1,...,x_k\}\in M$ there is a holomorphic map $g:{\Bbb C}\rightarrow M$, such that $g({\Bbb C})\supset\{x_1,...,x_k\}.$

Then $\widehat{s}_0(M)=\infty $.
\eT
To prove this we need the following two lemmas.
\bL\label{inf} Suppose $M$ is a complex manifold, $\dim (M)\geq 2$. Suppose also that 
for any distinct $k$ points $\{x_1,...,x_k\}\in M$ the following is true:

1. there is a holomorphic map $f:M\rightarrow {\Bbb C}$, such that $f(x_i)\neq
f(x_j)$ if $i\neq j$.

2. there is a holomorphic map $g:{\Bbb C}\rightarrow M$, such that $g({\Bbb C})\supset
\{x_1,...,x_k\}.$

Then $\widehat{s}_0(M)=\infty $.
\eL
\begin{pf} For any given $k$ points $\{x_1,...,x_k\}\in M$, fix $w_j\in
g^{-1}(x_j)$. Now consider $\Psi =g\circ \varphi \circ f:M\rightarrow M
$, where $\varphi :{\Bbb C}\rightarrow {\Bbb C}$ is the Lagrange polynomial, such that $
\varphi (f(x_j))=w_j$. Then $\Psi \neq id\,$ is a holomorphic endomorphism
of $M$ fixing the given points.
\end{pf}
\bL Let $M$ be a complex manifold such that for any two points $p\neq
q\in M$ there exists a holomorphic function $h:M\rightarrow {\Bbb C}$, such that $%
h(p)\neq h(q)$. Then for any finite number of distinct points $\{x_1,...,x_k\}\in M$ there exists a holomorphic function $f:M\rightarrow {\Bbb C}$
, that separates these points:  $f(x_i)\neq f(x_j)$ for $i\neq j$.
\eL
\begin{pf} Induction: suppose for $k\geq 2$ the statement holds. We note here that without any loss of generality we may assume that if a function separates $k$
given points, its values at these points can be preassigned as we please.
Now let the points $\{x_1,...,x_{k+1}\}\in M$ be given. For $m=1,...,k+1$ consider functions $f_m:M\rightarrow {\Bbb C}$ such that $f_m(x_s)=s$ for $s\neq m$. If no $f_m$ separates all $k+1$ points, then
for all $m$ the value  $f_m(x_m)\,$must be an integer (moreover  $f_m(x_m)\in \{1,...,k+1\}\backslash \{m\}$). Let $\alpha _1,...\alpha _{k+1}$
be a set of linearly independent numbers over ${\Bbb Z}$, the ring of integers. Consider $f=\sum\limits_{m=1}^{k+1}\alpha
_mf_m$. We claim that $f$ does the trick. Indeed for $i\neq j$, $
f(x_i)-f(x_j)=\sum\limits_{m\neq i,j}\alpha _m(i-j)+\alpha _i(f_i(x_i)-j)+\alpha
_j(i-f_j(x_j))\neq 0$ since the number of non-zero coefficients (equal to $(i-j)$) is at least $k-1\ge 1$.
\end{pf}
By the definition of a Stein manifold any two points can be separated by a holomorphic function. Therefore by the above lemma any finite number of points in such a manifold can be separated. The proof of Theorem \ref{sm} follows now from Lemma \ref{inf}.

\vspace{2pt}
\noindent{\it Remark.} The property described in the above theorem leads to a natural question that seems to be open: find the necessary and sufficient conditions for a complex manifold $M$ to have the following geometric property: any finite number of points can be connected by an analytic curve on $M$.


\begin{thebibliography}{999}

\bibitem[Bo]{Bo} N.~Bourbaki, Int$\grave{e}$gration, Hermann, Paris, 1963
\bibitem[Ca1]{Ca1} H.~Cartan, Les fonctions de deux variables complexeses et le probl$\grave{e}$me de la repr\'esentation analytique, {\it J. Math. pures et appl.}, $9^e$ s\'erie, 11 (1931) 1-114. 
\bibitem[Ca2]{Ca2} H.~Cartan, Sur les fonctions de plusieurs variables complexes. L'it$\acute{e}$ration des transformations int$\acute{e}$rieures d'un domaine born$\acute{e}$, {\it Math. Z. 35 (1932) 760-773.}

\bibitem[FF]{FF} S.~D.~Fisher and J.~Franks, The fixed points of an
analytic self-mapping, {\it Proc.\ AMS},  99(1987), 76--78.

\bibitem[FK1]{FK1}  B.~L.~Fridman, K.~T.~Kim, S.~G.~Krantz, \& D.~Ma, 
On fixed points and determining sets for holomorphic automorphisms, 
{\it Michigan Math.\ J.} 50(2002), 507--515. 


\bibitem[FK2]{FK2} B.~L.~Fridman, K.~T.~Kim, S.~G.~Krantz, \& D.~Ma, 
On Determining Sets for Holomorphic Automorphisms, to appear in
{\it Rocky Mountain J. of Math}. 

\bibitem[FM]{FM} B.~L.~Fridman, D.~Ma, Properties of Fixed Point Sets and a Characterization of the Ball in ${\Bbb C}^n$ , to appear in {\it Proc. AMS}.

\bibitem[GK]{GK} R.~E.~Greene and S.~G.~Krantz,
Characterization of complex manifolds by the isotropy subgroups of their automorphism groups, {\it Indiana Univ. Math. J. 34 (1985), no. 4, 865-879}.

\bibitem[GKM]{GKM} D. Gromoll, W. Klingenberg, and W. Meyer,
Riemannsche Geometrie im Grossen, $2^{\rm nd}$ ed., {\it Lecture Notes in
Mathematics}, v.~55,
Springer-Verlag, New York, 1975.

\bibitem[GR]{GR} R.~C.~Gunning, H.~Rossi, Analytic Functions of Several Complex Variables, Prentice-Hall, Inc., Englewood Cliffs, N.J. 1965

\bibitem[KK]{KK} K.~T.~Kim, S.~G.~Krantz, Determining Sets and Fixed Points for Holomorphic Endomorphisms, {\it Contemporary Math.} 328 (2003), 239-246.

\bibitem[K1]{K1} W.~Klingenberg, Riemannian Geometry,
$2^{\rm nd}$ ed., {\it de Gruyter Studies in Mathematics}, Berlin, 1995.

\bibitem[Ko]{Ko} S.~Kobayashi, Hyperbolic manifolds and holomorphic mappings, Marcel dekker, New York, 1970.


\bibitem[Le]{Le} K.~Leschinger, \"Uber fixpunkte holomorpher
Automorphismen, {\it Manuscripta Math.}, 25 (1978), 391-396.

\bibitem[Ma]{Ma} D.~Ma, Upper semicontinuity of Isotropy and
automorphism groups, {\it Math.\ Ann.}, 292(1992), 533--545.

\bibitem[Mas]{Mas} B. Maskit, The conformal group of a plane
domain,  {\it Amer.\ J.\ Math.}, 90 (1968), 718--722.


\bibitem[PL]{PL} E. Peschl and M. Lehtinen, A conformal self-map which
fixes 3 points is the identity, Ann.\ Acad.\ Sci.\ Fenn., Ser. A I Math.,  4
(1979), no.~1, 85--86.

\bibitem[Su]{Su} N. Suita, On fixed points of conformal
self-mappings, {\t Hokkaido Math.\ J.}, 10(1981), 667--671.

\bibitem[Vi1]{Vi1} J.-P. Vigu\'e, Sur les ensembles d'unicit\'e pour 
les automorphismes analytiques d'un domaine born\'e, {\it C.\ R.\ Acad.\ 
Sci.\ Paris, Ser.~I} 336(2003), 589--592. 

\bibitem[Vi2]{Vi2} J.-P. Vigu\'e, Ensembles d'unicit\'e pour 
les automorphismes et les endomorphismes analytiques d'un domaine born\'e, {\it Annales Institut Fourier},  55  (2005), p. 147-159.
\end{thebibliography}
\end{document}